# DISCUSSION: ONE-STEP SPARSE ESTIMATES IN NONCONCAVE PENALIZED LIKELIHOOD MODELS

By Peter Bühlmann and Lukas Meier

*ETH Zürich*

Hui Zou and Runze Li ought to be congratulated for their nice and interesting work which presents a variety of ideas and insights in statistical methodology, computing and asymptotics.

We agree with them that one- or even multi-step (or -stage) procedures are currently among the best for analyzing complex data-sets. The focus of our discussion is mainly on high-dimensional problems where $p \gg n$: we will illustrate, empirically and by describing some theory, that many of the ideas from the current paper are very useful for the $p \gg n$ setting as well.

**1. Nonconvex objective function and multi-step convex optimization.** The paper demonstrates a nice, and in a sense surprising, connection between difficult nonconvex optimization and computationally efficient Lasso-type methodology which involves one- (or multi-) step convex optimization. The SCAD-penalty function [5] has been often criticized from a computational point of view as it corresponds to a nonconvex objective function which is difficult to minimize; mainly in situations with many covariates, optimizing SCAD-penalized likelihood becomes an awkward task.

The usual way to optimize a SCAD-penalized likelihood is to use a local quadratic approximation. Zou and Li show here what happens if one uses a local *linear* approximation instead. In 2001, when Fan and Li [5] proposed the SCAD-penalty, it was probably easier to work with a quadratic approximation. Nowadays, and because of the contribution of the current paper, a local linear approximation seems as easy to use, thanks to the homotopy method [12] and the LARS algorithm [4]. While the latter is suited for linear models, more sophisticated algorithms have been proposed for generalized linear models; cf. [6, 8, 13].

In addition, and importantly, the local linear approximation yields sparse model fits where quite a few or even many of the coefficients in a linear or

---

generalized linear model are zero, that is, the method does variable selection. From this point of view, the local linear approximation is often to be preferred. In fact, it closely corresponds to the adaptive Lasso [17] which is, in our view, very useful for variable selection with Lasso-type technology. The rigorous convergence results in Section 2.3 of the paper, with a nice ascent property as for the EM-algorithm, are further reasons to favor the local linear approximation over the local quadratic versions with its heuristic rule for ad-hoc thresholding to zero [as described in the paragraph after formula (2.4) of the paper]. Finally, in Theorem 2 of the paper, the local linear approximation is shown, to yield the best convex majorization of the penalty function.

1.1. *Connection to the adaptive Lasso for type* 1 *penalty functions.* Section 4 in the paper distinguishes the cases where the regularization parameter can be separated from the penalty function or not. Type 1 penalty functions $p_\lambda(t) = \lambda p(t)$ allow for separation, for example, the Bridge penalties and the logarithm penalty, but excluding SCAD. From a computational point of view, the type 1 penalties are to be preferred because path-following algorithms can be used: this is Algorithm 1 in Section 4 of the paper, allowing to compute the whole regularization path very efficiently. In contrast, Algorithm 2, which can be used for the one-step SCAD estimator, seems much less efficient for approximating the entire regularization paths.

QUESTION. Why should we use the one-step SCAD estimator? In particular, (at least some of the) one-step type 1 penalty estimators, for example, the adaptive Lasso as discussed below, have the same asymptotic oracle properties, under the same conditions, as the one-step SCAD presented in Theorem 3.

Consider now type 1 penalty functions. Formulae (3.1) and (3.3) of the paper describe the one-step estimator based on the local linear approximation. In formula (3.1) corresponding to a linear model, we see that only the penalty function involves the initial estimator: the estimator can be written as

$$\beta^{(1)} = \underset{\beta}{\arg\min} \tfrac{1}{2}\|\mathbf{y} - \mathbf{X}\beta\|^2 + n\lambda \sum_{j=1}^{p} w(|\beta_j^{(0)}|)|\beta_j|,$$

where the penalty is based on re-weighting the $\ell_1$-norm (or Lasso-penalty) with weights $w_j = w(|\beta_j^{(0)}|)$ depending on the initial estimator. Of course, the weights also depend on the type 1 penalty function which we aim to approximate with the local linear approximation.

This is exactly the idea of the adaptive Lasso, recently proposed by [17]. We think it is important to emphasize this connection (which is not mentioned at all in the paper) because: (i) it is a simple and very effective idea to reduce the bias of the Lasso, see below; (ii) the adaptive Lasso is theoretically supported [17] and enjoys the same oracle result as the one-step SCAD described in Theorem 3, and there are theoretical results even for the high-dimensional situation where $p \gg n$ [7].

Regarding issue (i), particularly in cases with many ineffective (or non-substantial) covariates, the prediction optimal Lasso typically needs a large penalty parameter to get rid of these ineffective covariates. But a large penalty parameter implies substantial shrinkage to zero even for the coefficients corresponding to the important covariates. Solutions to improve such bias problems of the Lasso are based on two-stage procedures, for example, the LARS-OLS hybrid [4], the adaptive Lasso [17] or the relaxed Lasso [9].

## 2. Variable selection in the high-dimensional case.
For simplicity, consider a linear model

$$\mathbf{y} = \mathbf{X}\beta + \varepsilon, \tag{1}$$

as discussed in Section 3.1 of the paper. Our goal is variable selection (which is in many applications more relevant than prediction; we will present an example from biology in our Section 3.2). In high-dimensional problems where $p \gg n$, computational aspects become crucial. Since there are $2^p$ sub-models, we cannot inspect all of them (even when using efficient branch-and-bound methods). Ad-hoc methods can be used, but they may be very unstable, yielding poor results (e.g., forward variable selection); see [2]. On the other hand, having provably correct algorithms or methods, as the one in the paper with provable properties, is much more desirable.

The Lasso and its modifications belong to the latter class of methods. Regarding the computational feasibility, for linear models, the complexity to compute all sub-models from Lasso is $O(np\min(n,p))$ which is linear in the dimensionality $p$ if $p \gg n$. An important question is whether such computationally efficient estimators have good, provable statistical properties. Meinshausen and Bühlmann [10] showed consistency of the Lasso for variable selection in the high-dimensional setting where $p \gg n$: there is one major assumption, the neighborhood stability condition, which was shown to be sufficient and "almost" necessary (the wording "almost" refers to a numerical value which has to be $< 1$ for sufficiency and $\leq 1$ for necessity). Later, the irrepresentable condition has been worked out [16, 17] which is easier to interpret but is equivalent to the Meinshausen–Bühlmann assumption on neighborhood stability. The irrepresentable assumption is restrictive and easily fails to hold if the design matrix exhibits a too strong "degree of linear dependence" (of course, there is always linear dependence if $p \gg n$) or a too strong population absolute correlation among the covariates.

Since the irrepresentable (or neighborhood stability) condition is restrictive, one would like to understand Lasso's behavior under weaker assumptions: recently, consistency results in terms of

$$(2) \qquad \|\hat{\beta}(\lambda) - \beta\|_q = o_P(1) \qquad (n \to \infty), \ q \in \{1, 2\}$$

have been achieved; see [11, 14, 15]. The result in (2) has implications for variable selection. In case of fixed dimension $p$, (2) implies if $\beta_j \neq 0$, then $\hat{\beta}_j(\lambda) \neq 0$ with probability tending to 1 [because otherwise the convergence to zero in (2) would not hold]. That is,

(3)    The Lasso yields a too large model which contains the true model with high probability (tending to 1 as $n \to \infty$).

Under suitable conditions, the statement in (3) also holds in the high-dimensional case where $p \gg n$; see [11]. In addition, we point out that

(4)    The prediction optimal (w.r.t. MSE) tuned Lasso contains the true model with high probability (tending to 1 as $n \to \infty$).

This has been proved for simple cases in [10].

Putting the two facts (3) and (4) together, we can view the Lasso as an excellent and computationally efficient tool for "variable filtering," in the sense that the true model is with high probability a subset of the Lasso-estimated model. To appreciate the value of such a result, imagine that we have $p \approx 10'000$ and $n \approx 50$ (e.g., from microarray data). As the size of the Lasso-estimated sub-model is bounded by $\min(n, p)$, which equals $n$ if $p \gg n$, we pursue an immense dimensionality reduction from $p \approx 10'000$ to something of the order 50.

When viewing the Lasso as a variable filtering method, it is clear that we want to do an additional step (similar to the one-step estimator in the paper) which aims to go from the Lasso-estimated model in the first stage to the true model in a second stage. We have already touched upon two-stage procedures for addressing the bias problem in Lasso. The main proposals are the LARS-OLS hybrid [4], the relaxed Lasso [9] and the adaptive Lasso [17] with the Lasso as initial estimator: all of them reduce the bias and, in fact, this is what will lead to consistency in variable selection. We think (based on empirical evidence) that the latter, which is essentially the one-step estimator in the paper when using the Lasso as initial estimator, is a very elegant way to address Lasso's overestimation behavior. In addition, some theory for consistency in variable selection has been worked out for the high-dimensional case [7].

**3. Beyond the one-step estimator.** For regularization in high-dimensional spaces, we may want to use more than one or two regularization parameters. This can be achieved by pursuing more iterations where every iteration involves a separate tuning parameter (and as described below, those parameters are "algorithmically" constrained). We propose here the

Multi-Step Adaptive Lasso (MSA-LASSO):

1. Initialize the weights $w_j^{(0)} \equiv 1$ $(j = 1, \ldots, p)$.
2. For $k = 1, 2, \ldots, M$,
   use the adaptive Lasso with penalty function

   $$\lambda_*^{(k)} \sum_{j=1}^{p} w_j^{(k-1)} |\beta_j|,$$

   where $\lambda_*^{(k)}$ is the regularization parameter leading to prediction optimality. Denote the estimator by $\beta^{(k)} = \beta^{(k)}(\lambda_*^{(k)})$. In practice, the value $\lambda_*^{(k)}$ can be chosen via some cross-validation scheme.
   Up-date the weights

   $$w_j^{(k)} = \frac{1}{|\beta^{(k-1)}(\lambda_*^{(k-1)})_j|}, \qquad j = 1, \ldots, p.$$

For $k = 1$, we do an ordinary Lasso fit and $k = 2$ corresponds to the adaptive Lasso. Note that what is termed "one-step" in the paper corresponds here to $k = 2$. Note that Zou and Li initialize with $w_j^{(0)} \equiv 0$ (in the terminology of MSA-LASSO), yielding the MLE (in step $k = 1$). We find it more natural, and actually essential in the high-dimensional case with $p \gg n$, to initialize with the nonzero weights allowing for regularized fitting in step $k = 1$.

We will illustrate below the MSA-LASSO on a small simulated model and a real data set from molecular biology. Before, we describe some properties of the method which are straightforward to derive.

PROPERTY 1. MSA-LASSO increases the sparsity in every step in terms of the $\ell_0$-norm, that is, fewer selected variables in every step. As "heuristics," which is derived from the Zou and Li paper, MSA-LASSO is related to approximating the nonconvex optimization problem with the log-penalty $\sum_{j=1}^{p} \log(|\beta_j|)$; see the formula appearing just before Proposition 2.

PROPERTY 2. MSA-LASSO can be computed using the LARS algorithm for every step. The computational complexity of MSA-LASSO is bounded by $O(Mnp \min(n, p))$; due to the increase of sparsity, a later step is faster to compute than an early one. The computational load is in sharp contrast to computing all solutions over all $M$ steps when allowing for any $\lambda$ for each Lasso path: this would require many more essential operations.

MSA-LASSO is different from the multi-step procedure as analyzed in Section 2.3 of the paper: there, the regularization parameter $\lambda$ is fixed. We will also discuss below in Tables 1 and 2, for a simulated example, that the algorithmic restriction of choosing the regularization parameters in a sequentially optimal fashion seems very reasonable.

3.1. *Small simulation study.* To illustrate the proposed MSA-LASSO method, a small simulation study is carried out. We use a linear model as in (1) with covariates $X$ from a multivariate normal distribution with correlation matrix $\Sigma_{i,j} = \rho^{|i-j|}$ (for various values of $\rho$). The true underlying parameter vector is of the form $\beta = (c, \ldots, c, 0, \ldots, 0)^T$ with $p_{\text{act}}$ nonzero entries and $c$ such that the signal-to-noise ratio is 9 (which we find more relevant for practical applications than a signal-to-noise ratio of 21.25 as in example 1 in the paper). The number of predictors is set to $p = 500$. We choose the number of active variables $p_{\text{act}} \in \{3, 25\}$. In each simulation run, a training set of size 100 and a validation set of size 50 is used to determine the prediction optimal estimator. A total of 100 simulation runs are carried out for each parameter setting.

As performance measures we use the squared error $\|\hat{\beta} - \beta\|_2^2$ and the number of false positives (FP) $\sum_{j=1}^p I(\hat{\beta}_j \neq 0, \ \beta_j = 0)$. Table 1 illustrates the results for the case $p_{\text{act}} = 3$. We denote by 1-step $(k = 2)$ the MSA-LASSO with $k = 2$: it equals the adaptive Lasso with Lasso as initial estimator and we use the terminology "1-step" to be more consistent with the paper. Furthermore, 1-step opt. is the adaptive Lasso with Lasso as initial estimator (as for 1-step), but we optimize over a large 2-dimensional grid of the two regularization parameters which are involved in the initial Lasso and the adaptive Lasso step. Several conclusions can be made. The estimation error can be slightly decreased by an additional step if the correlation is not too high. More importantly, the number of false positives gets reduced. The number of false negatives is zero for this setting (not shown), that is, the true variables are always identified. The computational extra effort of the 1-Step opt. estimator does not pay off in this situation.

TABLE 1
*Results for $p_{\text{act}} = 3$ active variables. Average and standard deviations (in parentheses) of squared errors and of false positives (FP)*

|  | $\boldsymbol{\rho = 0}$ | $\boldsymbol{\rho = 0.2}$ | $\boldsymbol{\rho = 0.5}$ | $\boldsymbol{\rho = 0.8}$ |
|---|---|---|---|---|
|  | Squared error | | | |
| 1-Step (MSA-LASSO $k = 2$) | 0.12 (0.16) | 0.09 (0.09) | 0.11 (0.11) | 0.21 (0.23) |
| 2-Step (MSA-LASSO $k = 3$) | 0.08 (0.12) | 0.08 (0.10) | 0.09 (0.09) | 0.23 (0.30) |
| 1-Step opt. | 0.08 (0.09) | 0.09 (0.10) | 0.10 (0.09) | 0.19 (0.17) |
|  | False positives (FP) | | | |
| 1-Step (MSA-LASSO $k = 2$) | 6.10 (12.07) | 3.82 (6.59) | 3.73 (5.00) | 3.19 (6.23) |
| 2-Step (MSA-LASSO $k = 3$) | 2.93 (7.53) | 2.19 (6.24) | 1.65 (2.69) | 2.19 (5.09) |
| 1-Step opt. | 2.92 (7.13) | 3.08 (7.55) | 2.88 (6.00) | 3.36 (5.56) |

TABLE 2
*Results for 25 active variables. Average and standard deviations (in parentheses) of squared errors and of false positives (FP)*

|  | $\rho = 0$ | $\rho = 0.2$ | $\rho = 0.5$ | $\rho = 0.8$ |
|---|---|---|---|---|
|  | | Squared error | | |
| 1-Step (MSA-LASSO $k = 2$) | 6.15 (1.45) | 3.23 (1.06) | 1.78 (0.45) | 1.57 (0.32) |
| 2-Step (MSA-LASSO $k = 3$) | 6.25 (1.50) | 3.28 (1.12) | 1.92 (0.48) | 1.73 (0.35) |
| 1-Step opt. | 6.03 (1.49) | 2.95 (1.04) | 1.48 (0.49) | 1.25 (0.33) |
|  | | False positives (FP) | | |
| 1-Step (MSA-LASSO $k = 2$) | 32.91 (23.66) | 30.26 (18.16) | 17.91 (16.37) | 8.42 (9.01) |
| 2-Step (MSA-LASSO $k = 3$) | 27.43 (23.54) | 25.52 (17.44) | 14.60 (15.76) | 6.35 (6.73) |
| 1-Step opt. | 31.55 (22.11) | 24.55 (15.55) | 9.48 (11.37) | 2.76 (5.61) |

The results for $p_{\text{act}} = 25$ are given in Table 2. There is a slight loss in terms of mean squared error (MSE) for this setting when doing an additional step. Already, the 1-Step estimator loses compared to the initial Lasso estimator (not shown), in terms of MSE; likewise, the 1-Step opt. estimator has worse performance than the initial Lasso estimator. However, the number of false positives (FP) gets reduced again due to the increased sparsity. In terms of FP, the 1-Step opt. estimator seems to perform better for moderate and large values of $\rho$; but an additional step ($k = 4$) in MSA-LASSO would improve performance with respect to FP as well.

3.2. *Real data example from biology.* Reducing the number of false positives can be very desirable in biological applications since follow-up experiments can be costly and laborious. In our experience, it is often appropriate to do estimation on the conservative side with a low number of false positives because we still see more positives than what can be typically validated in a laboratory.

We illustrate the MSA-LASSO method on a problem of motif regression [3] for finding transcription factor binding sites in DNA sequences. [1] contains a collection of microarray data and a collection of motif candidates for yeast. The idea is to predict the gene expression value of a gene based on the corresponding motif scores (the information based on the sequence data). The data set which we consider consists of $n = 2587$ gene expression values of a heat-shock experiment and $p = 666$ motif scores. We use a training set of size 1300 and a validation set of size 650. The remaining data is used as a test-set.

The squared prediction error on the test-set $\mathbb{E}[(\hat{\mathbf{y}}_{\text{new}} - \mathbf{y}_{\text{new}})] = (\hat{\beta} - \beta)^T \Sigma (\hat{\beta} - \beta) + \text{Var}(\varepsilon)$ remains essentially constant for all estimators [probably due to high noise, that is, large value of $\text{Var}(\varepsilon)$]: 0.6193, 0.6230, 0.6226 for the Lasso, 1-Step and 2-Step estimator, respectively. But the number of selected variables decreases substantially:

|                             | Lasso | 1-Step | 2-Step |
|-----------------------------|-------|--------|--------|
| number of selected variables | 91    | 42     | 28     |

The list of top-candidate motifs gets slightly rearranged between the different estimators. The hope (and in part a verified fact) is that the 1- or 2-Step estimator yield more stable lists with fewer false positives.

Seminar für Statistik
ETH Zürich
CH-8092 Zürich
Switzerland
E-mail: buhlmann@stat.math.ethz.ch
           meier@stat.math.ethz.ch